\documentclass[11pt]{article}
\usepackage{amssymb,amsmath,graphicx,theorem}
\usepackage{hyperref}

\def\NN{{\cal N}}

\def\GG{{\cal G}}

\def\SS{{\cal S}}

\def\LL{{\cal L}}

\def\E{{\sf E}}

\def\PP{{\cal P}}

\def\R{{\mathbf R}}
\def\Z{{\mathbf Z}}

\def\T{{^{\sf T}}}

\def\eps{\varepsilon}

\def\diag{\hbox{\rm diag}}
\def\hom{\hbox{\rm hom}}
\def\inj{\hbox{\rm inj}}

\def\rk{\hbox{\rm rk}}

\def\simp{\text{\rm simp}}
\def\Eu{\text{\sf eul}}
\def\Pm{\text{\sf perf}}
\def\Tu{\text{\sf tut}}
\def\Cr{\text{\sf chr}}
\def\Fl{\text{\sf flo}}
\def\Ed{\text{\sf expt}}

\newtheorem{theorem}{Theorem}[section]
\newtheorem{prop}[theorem]{Proposition}
\newtheorem{lemma}[theorem]{Lemma}
\newtheorem{claim}{Claim}[section]

\newtheorem{corollary}[theorem]{Corollary}
\theorembodyfont{\rmfamily}

\long\def\killtext#1{}

\newenvironment{proof}{\noindent{\bf Proof. }}{\hfill$\square$\medskip}

\begin{document}

\title{Contractors and connectors of graph algebras\footnote{AMS
Subject Classification: Primary 05C99, Secondary 16S99}}
\author{{\sc L\'aszl\'o Lov\'asz} and {\sc Bal\'azs Szegedy}\\
Microsoft Research \\
One Microsoft Way}

\date{April 2005}

\maketitle

\tableofcontents

\begin{abstract}
We study generalizations of the ``contraction-deletion'' relation of
the Tutte polynomial, and other similar simple operations, to other
graph parameters. The question can be set in the framework of graph
algebras introduced by Freedman, Lov\'asz and Schrijver in
\cite{FLS}, and it relates to their behavior under basic graph
operations like contraction and subdivision.

Graph algebras were introduced in \cite{FLS} to study and
characterize homomorphism functions. We prove that for homomorphism
functions, these graph algebras have special elements called
``contractors'' and ``connectors''. This gives a new characterization
of homomorphism functions.
\end{abstract}

\section{Introduction and results}

The contraction-deletion operation for the Tutte polynomial is a
basic tool in graph theory. For our purposes, let us formulate this
as property as follows: Let $G$ be a graph and $u$ and $v$
nonadjacent nodes in $G$. Let $G'$ be obtained by identifying these
nodes. Then the Tutte polynomial of $G'$ can be expressed as a linear
combination of the Tutte polynomials of $G$ and the graph $G+uv$
(obtained by connecting $u$ and $v$ by an edge).

Which other graph parameters have a similar property that the
parameter of $G'$ can be expressed as a linear combination of the
parameter on graphs obtained from $G$ by attaching various ``small''
graphs at $u$ and $v$?

If we study the number of perfect matchings in a graph, then a useful
observation is that subdividing an edge by two new nodes does not
change this number. Which other graph parameters have a similar
property that the parameter of $G$ can be expressed as a linear
combination of the parameter on graphs obtained from $G$ by deleting
the edge $uv$ and attaching various ``small'' graphs at $u$ and $v$?

These questions are related to the work in \cite{FLS} and
subsequent work \cite{LL,LTS}. Here certain algebras generated by
graphs played a useful role, and the above questions can be stated
as rather basic properties of these algebras. Among others, they
can phrased in terms of the existence of special elements called
``contractors'' and ``connectors''.

Graph algebras were introduced in \cite{FLS} to study and
characterize homomorphism functions. We prove that for homomorphism
functions, these graph algebras have contractors and connectors. This
gives a new characterization of homomorphism functions.

\subsection{Graph algebras}

To state our results, we need to introduce some formalism. Fix a
positive integer $k$. A {\it $k$-labeled graph} is a finite graph in
which some of the nodes are labeled by numbers $1,\dots,k$. We denote
by $K_k$ the $k$-labeled complete graph with $k$ nodes, and by $O_k$,
the $k$-labeled graph with $k$ nodes and no edges. A $k$-labeled
quantum graph is {\it simple}, if it has no multiple edges, and its
labeled nodes are independent.

Let $F_1$ and $F_2$ be two $k$-labeled graphs. Their {\it product}
$F_1F_2$ is defined as follows: we take their disjoint union, and
then identify nodes with the same label. For two $\emptyset$-labeled
graphs, $F_1F_2$ is their disjoint union. Clearly this multiplication
is associative and commutative.

A {\it $k$-labeled quantum graph} is a formal finite linear
combination (with real coefficients) of $k$-labeled graphs. Let
$\GG_k$ denote the (infinite dimensional) vector space of all
$k$-labeled quantum graphs. We can turn $\GG_k$ into an algebra by
using $F_1F_2$ introduced above as the product of two generators,
and then extending this multiplication to the other elements
linearly. Clearly $\GG_k$ is associative and commutative. The
graph $O_k$ is a unit element in $\GG_k$.

We'll also consider the subalgebra $\GG^{\simp}_k$ generated by
$k$-labeled simple graphs and the subalgebra $\GG^0_k$ generated
by those $k$-labeled graphs whose labeled points are independent.
It is clear that $\GG^{\simp}_k\subseteq\GG^0_k\subseteq\GG_k$.

\subsection{Graph parameters and star algebras}

A {\it graph parameter} is a function defined on graphs, invariant
under isomorphism. We call a graph parameter $f$ {\it multiplicative}
if for any two (0-labeled) graphs $F_1,F_2$ we have
\[
f(F_1F_2)=f(F_1)f(F_2).
\]

Every graph parameter $f$ introduces further structure on the
algebras $\GG_k$. We extend $f$ linearly to quantum graphs, and
consider $f(x)$ as a ``trace'' of $x$.  We use this trace function to
introduce an inner product on $\GG$ by
\[
\langle x,y\rangle = f(xy).
\]
Let $\NN_k(f)$ denote the kernel of this inner product, i.e.,
\[
\NN_k(f)=\{x\in\GG_k:~f(xy)=0~\forall y\in\GG_k\}.
\]
Then we can define the factor algebra
\[
\GG_k/f=\GG_k/\NN_k(f).
\]
For $x,y\in\GG_k$, we write $x\equiv y\pmod f$ if $x-y\in\NN_k(f)$.

The dimension of $\GG_k/f$ is called the {\it rank-connectivity} of
the parameter $f$, and is denoted by $\rk(f,k)$. This is in general
infinite, but it is finite for quite a few interesting graph
parameters.

We say that $f$ is {\it reflection positive} if this inner product is
semidefinite: $\langle x,x\rangle\ge 0$ for all $x$. In this case, we
have
\[
\NN_k(f)=\{x\in\GG_k:~f(x^2)=0\}.
\]
Furthermore the algebra $\GG_k/f$ is a commutative algebra whose
elements form a Hilbert space with the property that $\langle
xy,z\rangle=\langle x, yz\rangle$. If this Hilbert space is finite
dimensional, then this implies that $\GG_k/f$ has a basis
$p_1,\dots,p_N$ such that $p_i^2=p_i$ and $p_ip_j=0$ for $i\not=
j$.

In this paper we prove a number of facts relating these graph
algebras to the basic graph operations of subdivision and
contraction. This will lead to alternate characterizations of
``homomorphism functions'' defined in Section \ref{HOMOMORPH}.

\subsection{Homomorphism functions}\label{HOMOMORPH}

This important class of graph parameters was the motivating example
for the studies in this paper. For two graphs $F$ and $G$, let
$\hom(F,G)$ denote the number of homomorphisms (adjacency preserving
maps) from $V(F)$ to $V(G)$.

We need to generalize this to the case when $G$ is weighted. A {\it
weighted graph} $G$ is a graph with a weight $\alpha_G(i)$ associated
with each node and a weight $\beta_G(i,j)$ associated with each edge
$ij$. In this paper we assume that the nodeweights are positive. Let
$\alpha_G=\sum_{i\in V(G)}$ denote the total nodeweight of $G$. (An
unweighted graph can be considered as a weighted graph where all the
node- and edgeweights are 1.)

To every function $\phi:~S\to V(G)$ with $S\subseteq V(F)$ we
assign the weight

\[
\hom_{\phi}(F,G)=\sum_{\psi:~V(F)\to
V(G)~,~\psi|_{S}=\phi}~~\prod_{i\in V(F)\setminus
S}\alpha_G(\phi(i)) \prod_{u,v\in V(F)} \beta_G(\psi(u),\psi(v)).
\]
We then define
\[
\hom(F,G)=\hom_{\emptyset}(F,G).
\]
Sometimes it is more convenient to use the ``homomorphism density''
\[
t(F,G)=\frac{\hom(F,G)}{\alpha_G^{|V(F)|}}.
\]
If $G$ is unweighted, this specializes to
\[
t(F,G)=\frac{\hom(F,G)}{|V(G)|^{|V(F)|}}.
\]
We'll also consider the number $\inj(F,G)$ of injective homomorphisms
of $F$ into $G$, and its normalized version
\[
t_0(F,G)=\frac{\inj(F,G)}{(|V(G)|)_{|V(F)|}}
\]
(where $(n)_k=n(n-1)\dots(n-k+1)$).

The following theorem was proved in \cite{FLS}:

\begin{theorem}\label{HOMCHAR}
A graph parameter $f$ can be represented in the form $f=\hom(.,H)$
for some finite weighted graph $H$ on at most $q$ nodes if and only
if it is reflection positive and $\rk(f,k)\le q^k$ for all $k\ge 0$.
\end{theorem}

An exact formula for $\rk(f, k)$ for homomorphism functions was
obtained in \cite{LL}. To state it, we need a definition. Two nodes
$i$ and $j$ in a weighted graph $H$ are {\it twins}, if
$\beta_H(i,k)=\beta_H(j,k)$ for every node $k\in V(H)$. Twin nodes
can be merged without changing the homomorphism functions $t(.,H)$
and $\hom(.,H)$.

\begin{theorem}\label{HOMDIM}
If $f=\hom(.,H)$, and $H$ has no twin nodes, then $\rk(f,k)$ is the
number of orbits of the automorphism group of $H$ on the ordered
$k$-tuples of nodes.
\end{theorem}

\subsection{Contractors and connectors}

For a 2-labeled graph $F$ in which the two labeled nodes are
nonadjacent, let $F'$ denote the graph obtained by identifying the
two labeled nodes. The map $F\mapsto F'$ maps 2-labeled graphs to
1-labeled graphs. We can extend it linearly to get an algebra
homomorphism $x\mapsto x'$ from $\GG^{0}_2$ into $\GG_1$.

The map $x\mapsto x'$ does not in general preserve the inner
product or even its kernel; we say that the graph parameter $f$ is
{\it contractible}, if for every $x\in\GG^{0}_2$, $x\equiv 0\pmod
f$ implies $x'\equiv 0\pmod f$; in other words, $x\mapsto x'$
factors to a linear map $\GG_2^{0}/f\to\GG_1/f$.

We say that $z\in\GG_2$ is a {\it contractor for $f$} if for every
$x\in\GG^{0}_2$, we have
\[
f(xz)=f(x').
\]
Informally, attaching $z$ at two nodes acts like identifying those
two nodes.

Our second concern is to get rid of multiple edges. We say that
$z\in\GG^{\simp}_2$ is a {\it connector} for $f$, if $z\equiv
K_2\pmod f$, i.e., for every $x\in\GG_2$ we have
\[
f(zx)=f(K_2x).
\]
Note that $K_2$ is always a connector, but it is not simple in the
sense defined above.

If a graph parameter $f$ has a simple connector $z$, then every
$k$-labeled quantum graph is congruent to a simple quantum graph
modulo $f$. Indeed, for every $x\in\GG_k$, in every $k$-labeled graph
in the expansion of $x$, every edge can be replaced by the simple
connector, which creates a simple quantum graph. In other words,
$\GG_k^\simp/f=\GG_k/f$.

Several general facts about connectors and contractors will be
stated and proved in Section \ref{GENERAL}.

\subsection{The algebra of concatenations}

For two $2$-labeled graphs $F_1$ and $F_2$, we define their {\it
concatenation} by identifying node $2$ of $F_1$ with node $1$ of
$F_2$, and unlabeling this new node. We denote the resulting
$2$-labeled graph by $F_1\circ F_2$. It is easy to check that this
operation is associative (but not commutative). We extend this
operation linearly over $\GG_2$.

This algebra has a $*$ operation: for a $2$-labeled graph $F$, we
define  $F^*$ by interchanging the two labels. Clearly
$(F_1F_2)^*=F_2^*F_1^*$. We can also extend this linearly over
$\GG_2$.

Let $f$ be a graph parameter. It is easy to see that if $x\equiv
0\pmod f$ then $x^*\equiv 0\pmod f$, so the $*$ operator is well
defined on elements of $\GG_2/f$. An further important property of
concatenation is that
\[
f((x\circ y)z)=f(x(z\circ y^*))
\]
for any three elements $x,y,z\in\GG_2$. It follows that if $x\equiv
0~\pmod f$ then $x\circ y\equiv 0\pmod f$ for every $y\in\GG_2$ and
thus concatenation is also well defined on the elements of $\GG_2/f$.
It is easy to see that $\mathcal{A}_1=(\GG_2/f,+,\circ)$ is an
associative (but not necessarily commutative) algebra. Note that if
$x,y\in\GG_2$ then $x\circ y\in\GG^0_2$. It follows in particular
that if $\mathcal{A}_1$ has a unit element then $\GG^0_2/f =
\GG_2/f$.

\begin{lemma} If an element $z\in\GG_2$ is a contractor for $f$ then the image of
$z$ under the map $\GG_2\to \GG_2/f$ is the unit element of the
algebra $\mathcal{A}_1$.
\end{lemma}

\begin{proof}
We have to check that $z\circ x\equiv x~\pmod f$ for all $x\in\GG_2$.
This is equivalent with $f((z\circ x)y)=f(xy)$ for all $x,y\in\GG_2$.
Using that $x\circ y\in \GG^{0}_2$ we obtain that
\[
f((z\circ x)y)=f(z(y\circ x^*))=f((y\circ x^*)')=f(xy).
\]
\end{proof}

\subsection{Contractors and connectors for homomorphism functions}

The first two results in this paper concerns graph parameters that
are homomorphism functions, i.e., they are of the form $f=\hom(.,H)$
for some weighted graph $H$. It is easy to check from the definitions
that a 2-labeled quantum graph $z$ is a contractor for $\hom(.,G)$ if
and only if
\begin{equation}\label{CONT-HOMPHI}
\hom_\phi(z,G)=\begin{cases}
    1 & \text{if $\phi(1)=\phi(2)$}, \\
    0 & \text{otherwise,}
  \end{cases}
\end{equation}
for every $\phi:~\{1,2\}\to V(G)$. It is a connector for $\hom(.,G)$
if and only if $z\in\GG_2^\simp$, and
\begin{equation}\label{CONN-HOMPHI}
\hom_\phi(z,G)=\beta_G(\phi(1),\phi(2))
\end{equation}
for every $\phi:~\{1,2\}\to V(G)$.

We denote by $P_n$ the path with $n$ nodes, with the two endnodes
labeled $1$ and $2$ (so $P_2=K_2$). A {\it quantum path} is a linear
combination of such paths. A {\it series-parallel graph} is a
$2$-labeled graph obtained from $K_2$ by repeated application of the
product and concatenation operations. A {\it series-parallel quantum
graph} is a linear combination of series-parallel graphs.

\begin{theorem}\label{CONT-CONN}
Let $f=\hom(.,H)$ for some finite weighted graph $H$. Then $f$ has a
contractor and also a simple connector. Furthermore, it has a
contractor that is a series-parallel quantum graph and a simple
connector that is a quantum path.
\end{theorem}

Using the notion of a contractor, we can give the following
characterization of homomorphism functions (we don't know whether a
similar theorem holds using some special connectors instead of
contractors).

\begin{theorem}\label{CONTRACT}
A graph parameter $f$ can be represented in the form $f=\hom(.,H)$
for some finite weighted graph $H$ if and only if it is
multiplicative, reflection positive and has a contractor.
\end{theorem}

\subsection{Homomorphisms into measure graphs}

Every symmetric measurable function $W:~[0,1]^2\to[0,1]$ defines a
graph parameter $t(.,W)$ by
\[
t(F,W)=\int\limits_{[0,1]^n} \prod_{ij\in E(F)}
W(x_i,x_j)\,dx_1\,\dots\,dx_n.
\]
Homomorphism density functions into finite weighted graphs is a
special case. Call a symmetric function $W:~[0,1]^2\to[0,1]$ a {\it
step function}, if there is a partition $[0,1]=A_1\cup\dots\cup A_q$
into measurable sets such that $W$ is constant on $A_i\times A_j$ for
all $1\le i,j\le q$. It is trivial to check that if $W$ is a step
function, then $t(.,W)=t(.,H)$ for a finite weighted graph $H$ and
vice versa. It was noted in \cite{FLS} that the graph parameter
$t(.,W)$ is reflection positive, and it is obvious that it is
multiplicative.

These parameters occur in the context of limits of graph sequences.
We say that a sequence $(G_n)$ of simple graphs is {\it convergent}
if $t(F,G_n)$ converges to some value $t(F)$ for every simple graph
$F$. In \cite{LSz} it was shown that (at least for the case of
parameters defined on simple graphs), the parameters $t(.,W)$ are
precisely the limits of parameters $t(F)$ obtained this way.

\begin{theorem}\label{W-PROP}
The graph parameter $t(.,W)$ is contractible, but has no contractor
unless $W$ is a step function.
\end{theorem}

\section{Examples}\label{EXAMPLES}

The following examples are described in more detail in \cite{FLS}.
Here we only discuss those properties of them that relate to
contractors and connectors.

\subsection{Matchings}\label{MATCHING}

Let $\Pm(G)$ denote the number of perfect matchings in the graph $G$.
It is trivial that $\Pm(.)$ is multiplicative. Its
node-rank-connectivity is exponentially bounded,
\[
\rk(\Pm,k)=2^k,
\]
but $\Pm$ is not reflection-positive. Thus $\Pm(G)$ cannot be
represented as a homomorphism function.

On the other hand: $\Pm$ has a contractor: a path of length 2, and
also a simple connector: a path $P_4$ of length $3$.

\subsection{Chromatic polynomial}\label{CHROMPOLY}

Let $\Cr(G)=\Cr(G,x)$ denote the chromatic polynomial of the graph
$G$. For every fixed $x$, this is a multiplicative graph parameter.
For $k,q\in\Z_+$, let $B_{kq}$ denote the number of partitions of a
$k$-element set into at most $q$ parts. So $B_k=B_{kk}$ is the $k$-th
Bell number. With this notation, we have \cite{FLW}
\[
\rk(\Cr,k)=
  \begin{cases}
    B_{kx} & \text{if $x$ is a nonnegative integer}, \\
    B_k & \text{otherwise}.
  \end{cases}
\]
Note that this is always finite, but if $x\notin \Z_+$, then it grows
faster than $c^k$ for every $c$. Furthermore, $M(\Cr,k)$ is positive
semidefinite if and only if either $x$ is a positive integer or $k\le
x+1$. The parameter $M(\Cr,k)$ is reflection positive if and only if
this holds for every $k$, i.e., if and only if $x$ is a nonnegative
integer, in which case indeed $\Cr(G,x)=\hom(G,K_x)$.

This parameter has a contractor for every $x$: the 2-labeled quantum
graph $K_2-O_2$ (which amounts to the standard contraction-deletion
identity for the chromatic polynomial). It is not hard to check that
$\frac{1}{x-1}P_3-\frac{x-2}{x-1}P_2$ is a simple connector if
$x\not=1$; for $x=1$, the chromatic polynomial is $0$ if there is an
edge, so $P_2$ is a simple connector.

\subsection{Flows}\label{FLOWS}

Let $\Gamma$ be a finite abelian group and let $S\subseteq\Gamma$ be
such that $S$ is closed under inversion. For any graph $G$, fix an
orientation of the edges. An {\em $S$-flow} is an assignment of an
element of $S$ to each edge such that for each node $v$, the product
of elements assigned to edges entering $v$ is the same as the product
of elements assigned to the edges leaving $v$. Let $\Fl(G)$ be the
number of $S$-flows.  This number is independent of the orientation.
In the case when $S=\Gamma\setminus \{0\}$, $\Fl(G)$ is the number of
nowhere-0 $\Gamma$-flows.

The parameter $\Fl(G)$ can be described as a homomorphism function
\cite{FLS}. It has a trivial simple connector, a path of length 2
(which is an algebraic way of saying that if we subdivide an edge,
then the flows don't change essentially). In the case of nowhere-0
flows, $K_2+O_2$ is a contractor (which amounts to the
contraction-deletion identity for the flow polynomial), but in
general, there does not seem to be a simple explicit construction for
a contractor.

\subsection{Tutte polynomial}

Consider the following version of the Tutte polynomial: in terms of
the variables $q$ and $v$, we have
\[
\Tu(G;q,v)=\sum_{A\subseteq E(G)} q^{c(A)}v^{|A|},
\]
where $c(A)$ denotes the number of components of the graph
$(V(G),A)$. This differs from the usual Tutte polynomial $T(x,y)$
on two counts: first, instead of the standard variables $x$ and
$y$, we use $q=(x-1)(y-1)$ and $v=y-1$; second, we scale by
$q^{c(E)}v^{n-c(E)}$. This way we lose the covariance under
matroid duality; but we gain that the contraction/deletion
relation holds for all edges $e$:
\begin{equation}\label{TUTTESKEIN}
\Tu(G)=v\Tu(G/e)+\Tu(G\setminus e).
\end{equation}
If $i$ is an isolated node of $G$, then we have
\begin{equation}\label{TUTTEISO}
\Tu(G-i)=q\Tu(G).
\end{equation}
If $G$ is the empty graph (no nodes, no edges), then $\Tu(G)=1$.
Another way of expressing (\ref{TUTTESKEIN}) is that $K_2 - vO_2$ is
a contractor of $\Tu$. It is not hard to check that
$(1/v)P_3-(1+q/v)O_2$ is a simple connector.

The chromatic polynomial and the number of nowhere-$0$ $\Gamma$-flows
are special substitutions into the Tutte polynomial. More precisely,
\[
\Cr(G;x)=\Tu(G;x,-1),
\]
and the number of nowhere-0 $k$-flows is
\[
\Fl(G)=\frac{(-1)^{|E(G)|}}{k^{|V(G)|}}\Tu(G;k,-k).
\]
It can be shown \cite{FLW} that for $v\not=0$, the Tutte polynomial
behaves exactly as the corresponding chromatic
polynomial:
\[
\rk(\Cr,k)=
  \begin{cases}
    B_{kq} & \text{if $q$ is a nonnegative integer}, \\
    B_k & \text{otherwise}.
  \end{cases}
\]
Furthermore, $\Tu(G;q,v)$ is reflection positive if and only if $q$
is a positive integer. Theorem \ref{HOMCHAR} implies that in this
case $\Tu(G;q,v)$ is a homomorphism function, while for other
substitutions it is not.

\subsection{The role of multiple edges}\label{MULTIPLE}

Let, for each (multi)graph $G$, $\tilde G$ denote the (simple) graph
obtained from $G$ by keeping only one copy of each parallel class of
edges. Consider a random graph $H$ on $N$ nodes with edge probability
$1/2$, then
\[
\Ed(G)=\E(t_0(G,H))=2^{-|E(\tilde G)|}
\]
is independent of $N$. It is not hard to see that we also have with
probability 1
\[
\Ed(G)=\lim_{N\to\infty} t(G,H).
\]
From this (or from direct computation) it follows that this graph
parameter is multiplicative and reflection positive. It can be
checked that $\rk(\Ed,k)=2^{k\choose 2}$, which is finite for every
$k$, but has superexponential growth.

The graph parameter $\Ed$ is not contractible. Consider the 3-star
$S_4$ with 2 endnodes labeled and the path $P_4$ with 3 edges with
both endnodes labeled. Then $S_4\equiv P_4 \pmod f$, but identifying
the labeled nodes produces a pair of parallel edges in $S_4$ but not
in $P_4$, so $f(S_4')=1/4$ but $f(P_4')=1/8$, showing that
$S_4'\not\equiv P_4' \pmod f$. This implies by lemma \ref{CR2CE}
below that $\Ed$ does not have a contractor. It is easy to see that
$f$ does not have a simple connector either.

\subsection{The number of eulerian orientations}

Let $\Eu(G)$ denote the number of eulerian orientations of the graph
$G$. It was remarked in \cite{LSz} that this parameter can be
expressed as $\Eu(G)=t(G,W)$, where
\[
W(x,y)=2\cos(2\pi(x-y)).
\]
Thus it follows by Theorem \ref{W-PROP} that $\Eu$ is contractible,
but has no contractor. It is easy to see that a path of length 2 is a
simple connector.

\section{General facts about connectors and contractors}\label{GENERAL}

We start with an easy observation.

\begin{prop}\label{CR2CE}
If a graph parameter has a contractor, then it is contractible.
\end{prop}

\begin{proof}
Let $z$ be a contractor for $f$. Suppose that $x\in\GG_2$ satisfies
$x\equiv 0\pmod f$, and let $y\in\GG_1$. Choose a $\hat y\in\GG_2$
such that $\hat y'=y$. Then
\[
f(x'y)=f(x'\hat y')=f((x\hat y)')=f((x\hat y)z)=f(x(\hat y z))=0,
\]
showing that $x'\equiv 0\pmod f$.
\end{proof}

While the existence of a contractor does not imply the existence of a
simple connector or vice versa, there is some connection, as
expressed in the following proposition (see also Corollary
\ref{SDCONT}).

\begin{prop}\label{CE2CR}
If $f$ is contractible, has a simple connector, and $\rk(f,2)$ is
finite, then $f$ has a contractor.
\end{prop}

\begin{proof}
Since $\langle x,y\rangle=f(xy)$ is a symmetric (possibly
indefinite) bilinear form that is not singular on $\GG_2/f$, there
is a basis $p_1,\dots,p_r$ in $\GG_2/f$ such that $f(p_ip_j)=0$ if
$i\not= j$ and $f(p_ip_i)\not=0$. By the assumption that $f$ is
simplifiable, we may represent this basis by simple quantum
graphs; then the contracted quantum graphs $p_i'$ are defined. Let
\[
z=\sum_{i=1}^N \frac{f(p_i')}{f(p_i^2)}p_i.
\]
We claim that $z$ is a contractor. Indeed, let $x\in\GG^{\simp}_2$,
and write
\[
x\equiv\sum_{i=1}^N a_ip_i\pmod f.
\]
Then we have
\[
f(xz)=\sum_{i=1}^N a_i \frac{f(p_i')}{f(p_i^2)}f(p_i^2)=\sum_{i=1}^N
f(p_i')a_i.
\]
On the other hand, contractibility implies that
\[
x'\equiv\sum_{i=1}^N a_ip_i'\pmod f,
\]
and so
\[
f(x')=\sum_{i=1}^N a_if(p_i')=f(xz).
\]
\end{proof}

\begin{prop}\label{SUBDIV}
If $M(f,2)$ is positive semidefinite, $f$ is contractible and
$\rk(f,2)$ is finite, then $f$ has a simple connector that is a
quantum path.
\end{prop}

For the proof, we need the following simple lemma.

\begin{lemma}\label{SHORTEN}
Assume that $M(f,2)$ is positive semidefinite. Let $x\in\GG_2$,
and assume that $x\circ P_3\equiv 0\pmod{f}$. Then $x\circ
P_2\equiv 0\pmod{f}$.
\end{lemma}

\begin{proof}
We have
\[
(x\circ P_2)^2=(x\circ P_3)x=0,
\]
and by reflection positivity, this implies that $x\circ P_2\equiv
0\pmod{f}$.
\end{proof}

\noindent{\bf Proof} [of Proposition \ref{SUBDIV}]: Since
$\GG_2/f$ is finite dimensional, there is a linear dependence
between $P_2,P_3,\dots$ in $\GG_2/f$. Hence there is a (smallest)
$k\ge 2$ such that $P_k$ can be expressed as
\begin{equation}\label{PEXPR}
P_k\equiv \sum_{i=1}^N a_iP_{k+i}\quad\pmod f
\end{equation}
with some positive integer $N$ and real numbers $a_1,\dots,a_N$.
The assertion is equivalent to saying that $k=2$.

Let $x=P_2-\sum_{i=1}^N a_iP_{2+i}$. Then (\ref{PEXPR}) can be
written as $x\circ P_{k-1} \equiv 0 ~ \pmod{f}$. If $k>3$, then
Lemma \ref{SHORTEN} implies that $x\circ P_{k-2} \equiv 0 ~
\pmod{f}$, which contradicts the minimality of $k$. Suppose that
$k=3$. Then from (\ref{PEXPR}) we have that $(x-\sum_{i=1}^N
a_ix\circ P_{1+i})\circ P_3\equiv x\circ x\circ P_2\equiv 0
~\pmod{f}$. By Lemma \ref{SHORTEN} we get that $x\circ x\equiv
0~\pmod{f}$ and using contractibility we obtain that $0=f((x\circ
x)')=f(x^2)$. Now reflection positivity shows that $x\equiv 0
\pmod{f}$.

\begin{corollary}\label{SIMPLE}
If $M(f,2)$ is positive semidefinite, f is contractible and
$\rk(f,2)$ is finite, then $\GG_k=\GG_k^\simp$ for every $k\ge 1$.
\end{corollary}

\noindent The following statement is a corollary of Proposition
\ref{SUBDIV} and Proposition \ref{CE2CR}.

\begin{corollary}\label{SDCONT}
If $M(f,2)$ is positive semidefinite, $f$ is contractible and
$\rk(f,2)$ is finite then $f$ has a contractor.
\end{corollary}

\section{Homomorphism functions: proofs}

\subsection{Proof of Theorem \ref{CONT-CONN}}

Suppose that $f = \hom(.,G)$ for some weighted graph $G$. We may
assume that $G$ is twin-free.

We start with constructing a connector. Let $\alpha_1,\dots,\alpha_m$
be the nodeweights and $\beta_{ij}$ ($i,j=1,\dots,m$), the
edgeweights of $G$. Let $B=(\beta_{ij})$ be the (weighted) adjacency
matrix of $G$, and let
$D=\diag(\sqrt{\alpha_1},\dots,\sqrt{\alpha_m})$. Let
$\lambda_1,\dots,\lambda_t$ be the nonzero eigenvalues of the matrix
$DBD$ (which are real as $DBD$ is symmetric), and consider the
polynomial $\rho(z)=z\prod_{i=1}^t (1-z/\lambda_i)$. Then
$\rho(DBD)=0$. Since the constant term in $\rho(z)$ is $0$ and the
linear term is $z$, this expresses $DBD$ as a linear combination of
higher powers of $DBD$:
\[
DBD=\sum_{s=2}^t a_s(DBD)^s,
\]
or
\begin{equation}\label{DBD}
B=\sum_{s=2}^t a_s (BD^2)^{s-1}B.
\end{equation}

For every mapping $\varphi:~\{1,2\}\to V(G)$, we have
\[
\hom_\varphi(P_s,G)=((BD^2)^{s-2}B)_{\varphi(1)\varphi(2)}.
\]
Let
\[
y=\sum_{s=2}^t a_s P_{s+1},
\]
Then (\ref{DBD}) implies that for every 2-labeled graph $G$,
\[
f(K_2G)= f(yG).
\]
Thus $y$ is a connector. By construction, it is a linear combination
of paths.

For the existence of a contractor, there are two general arguments.

First, we can use Lemma \ref{CE2CR}: it is easy to check that $f$ is
contractible; the condition that $\GG_1/f=\GG_2'/f$ follows from the
existence of a connector; and $M(f,2)$ is positive semidefinite and
has finite finite rank by Theorem \ref{HOMCHAR}.

Second, to prove that there exists a 2-labeled quantum graph $z$
satisfying (\ref{CONT-HOMPHI}), we can invoke the following result
\cite{LL}:

\begin{theorem}
{\it Let $G$ be a twin-free weighted graph and $\Phi:~V(G)^k\to \R$.
Then there exists a $k$-labeled quantum graph $z$ such that
\[
\hom_\phi(z,G)=\Phi(\phi)
\]
for every  $\phi\in V(G)^k$, if and only if $\Phi$ is invariant under
the automorphisms of $G$: for every $\phi\in V(H)^k$ and every
automorphism $\sigma$ of $H$, $\Phi(\phi\sigma) =\Phi(\phi)$.}
\end{theorem}

However, it is worth while to give a third, more specific argument,
because it gives the stronger result that a series-parallel
contractor exists. To every $2$-labeled quantum graph $x$, we assign
the $V(G)\times V(G)$ matrix $M(x)$ as follows: for $i,j\in V(G)$, we
set
\[
M(x)_{ij}=\hom_{1\mapsto i,2\mapsto j}(x,G).
\]
Then it is easy to check that $M$ is a linear map from $\GG_2$ to the
space $V(G)\times V(G)$ matrices, and it also respects products in
the following sense:
\[
M(x\circ y)=M(x)D M(y), \qquad M(xy)=M(x)\circ M(y)
\]
(here $M(x)\circ M(y)$ denotes the Schur, or elementwise, product of
these matrices). Furthermore, interchanging the labels $1$ and $2$
corresponds to transposition of the corresponding matrix. Clearly,
$M(K_2)=B$ is the weighted adjacency matrix of $G$.

Now let $\SS\PP\subseteq \GG_2$ denote the space of series-parallel
quantum graphs, and let $\LL$ be the set of corresponding matrices.
We want to show that the identity matrix $I$ is in $\LL$. Clearly
$\LL$ is a linear space that is also closed under the Schur product,
the operation $(X,Y)\mapsto XDY$, and transposition. So the theorem
follows if we prove the following algebraic fact.

\begin{lemma}\label{ALGEBRA}
Let $\LL$ be a linear space of $n\times n$ matrices, and let $D$ be a
diagonal matrix with positive entries in the diagonal. Assume that
$\LL$ is closed under transposition, Schur product and the operation
$(X,Y)\mapsto XDY$. Assume furthermore that no row is 0 in every
matrix in $\LL$, and no two rows are parallel in every matrix in
$\LL$. Then $\LL$ contains the identity matrix.
\end{lemma}

\begin{proof}
We start with a remark.

\begin{claim}\label{FUNCTION}
Let $f:~\R\to\R$ be any function and $M\in\LL$. Then the matrix
$f(M)$, obtained by applying $f$ to every entry of $M$, is also in
$\LL$.
\end{claim}

Indeed, on the finite number of real numbers occurring as entries of
$M$, the function $f$ equals to some polynomial $\sum_{i=0}^N
a_ix^i$. Then
\[
f(M)=\sum_{i=0}^N a_i M^{(i)},
\]
where $M^{(i)}$ is the Schur product of $i$ copies of $M$. This shows
that $f(M)\in\LL$.

Clearly, there is a ``generic'' element $W\in\LL$ such that no row or
column of $W$ is $0$ and no two rows or columns of $W$ are equal.
Replacing $W$ by $W\circ W+\eps W$ with a small enough $\eps$, we may
also assume that $W\ge 0$.

We claim that all maximal entries of $W^TDW$ are on the diagonal.
Indeed, suppose that $(W^TDW)_{ij}$ is a maximal entry. By
Cauchy--Schwartz, we have
\begin{align*}
(W^TDW)_{ij}&=\sum_{h\in V(G)}\alpha_h W_{hi} W_{hj}\le
\Bigl(\sum_{h\in V(G)}\alpha_hW_{hi}^2\Bigr)^{1/2}\Bigl(\sum_{h\in
V(G)}\alpha_hW_{hj}^2\Bigr)^{1/2}\\
&=((W^TDW)_{ii})^{1/2}((W^TDW)_{jj})^{1/2}.
\end{align*}
It follows that $(W^TDW)_{ij}=(W^TDW)_{ii}=(W^TDW)_{jj}$, and that
the $i$-th column of $W$ is parallel to the $j$-th. Since $W\ge 0$,
this implies that the $i$-th column is equal to the $j$-th, and hence
by the choice of $W$ it follows that $i=j$.

Applying Claim \ref{FUNCTION}, we can replace the maximal entries
of $W\T DW$ by 1 and all the other entries by 0, to get a nonzero
diagonal 0-1 matrix $P\in\LL$. Choose such a matrix $P$ with
maximum rank; we claim that it is the identity matrix.

Suppose not, and consider the matrix $Q=I-P$ (we don't know yet
that $Q\in\LL$). The matrix $PDP=PD=DP$ is in $\LL$, and applying
Claim \ref{FUNCTION} again, we get that
$PD^{-1}P=PD^{-1}=D^{-1}P\in\LL$. Hence for every matrix
$M\in\LL$,
\begin{align*}
QMQ&=M-MP-PM+PMP\\
&=M-MD(D^{-1}P)-(PD^{-1})DM +(PD^{-1})DMD(D^{-1}P)\in\LL.
\end{align*}
In particular, $QW^TDWQ\in\LL$. By the same argument as above, we
see that all maximal entries of $QW^TDWQ$ are on its diagonal, and
so applying Claim \ref{FUNCTION} again, we get that $M$ contains a
nonzero matrix $Q'$ obtained from $Q$ by changing some of its 1's
to 0. Now $P+Q'\in\LL$ is a diagonal 0-1 matrix with larger rank
than $P$, a contradiction.
\end{proof}

\subsection{Proof of Theorem \ref{CONTRACT}}

The necessity of the conditions follows by Theorem \ref{CONT-CONN}.

To prove the sufficiency of the conditions, it suffices to prove that
there exists a $q>0$ such that $\rk(M(f,k))\le q^k$ for all $k\ge 0$,
and then invoke Theorem \ref{HOMCHAR}. Note that reflection
positivity is used twice: the existence of a contractor does not in
itself imply an exponential bound on the rank connectivity (cf.
Example \ref{CHROMPOLY}). Let $g_0$ be a contractor for $f$; we show
that $q=f(g_0^2)$ satisfies these conditions. Since $f$ is
multiplicative, we already know this for $k=0$.

We may normalize $f$ so that $f(K_1) = 1$. Let $N=\rk(M(f,k))$.
Consider the basic idempotents $p_1,\dots, p_N$ in the algebra of
$k$-labeled quantum graphs defined by $f$, and let $q_i =
\frac{1}{\sqrt{f(p_i)}}p_i$. Let $q_i \otimes q_i$ denote the
$(2k)$-labeled quantum graph obtained from $2k$ labeled nodes by
attaching a copy of $q_i$ at $\{1,\dots,k\}$ and another copy of
$q_i$ at $[k+1,2k]$. Let $h$ denote the $(2k)$-labeled quantum graph
obtained from $2k$ labeled nodes by attaching a copy of $g_0$ at
$\{i,k+i\}$ for each $i = 1,\dots,k$. Consider the quantum graph
\[
x = \sum_{i=1}^N q_i \otimes q_i - h.
\]
By reflection positivity, we have $f(x^2)\ge 0$. But
\[
f(x^2) = \sum_{i=1}^N \sum_{i=1}^N \langle q_i \otimes q_i, q_i
\otimes q_i \rangle -2 \sum_{i=1}^N \langle q_i \otimes q_i, h\rangle
+ \langle h,h\rangle.
\]
Here by the fact that $q_iq_j = 0$ if $i\not= j$, we have
\[
\sum_{i=1}^N \sum_{j=1}^N \langle q_i \otimes q_i, q_j \otimes
q_j\rangle = \sum_{i=1}^N \langle q_i , q_i\rangle^2 = N.
\]
Furthermore, by the definition of $g_0$ and $h$, we have
\[
\sum_{i=1}^N \langle q_i \otimes q_i, h\rangle= \sum_{i=1}^N \langle
q_i, q_i\rangle = N.
\]
Finally, by the definition of $h$ and the multiplicativity of $f$, we
have
\[
\langle h,h \rangle= f(g_0^2)^k.
\]
Thus $f(x^2)\ge 0$ implies that $N \le f(g_0^2)^k$, which completes
the proof.

\subsection{Proof of Theorem \ref{W-PROP}}

Set $t=t(.,W)$. Let $F:~[0,1]^2\to[0,1]$ be any integrable function.
We define
\[
\|F\|_x=\int\limits_0^1 |F(x,y)|\,dy.
\]
Then
\begin{equation}\label{L1X}
\|F\|=\int\limits_0^1 \|F\|_x \,dx
\end{equation}
is the usual $\ell_1$-norm of $F$.

\begin{lemma}\label{W-U}
Let $U,W:~[0,1]^2\to[0,1]$ be two symmetric functions and let $F$ be
a 2-labeled graph with $m$ edges in which the labeled nodes are
independent. Let $d_i$ denote the degree of $i$ in $F$. Then for
every $x,y\in[0,1]$,
\[
|t_{xy}(F,U)-t_{xy}(F,W)|\le d_1\|U-W\|_x +
d_2\|U-W\|_y+(m-d_1-d_2)\|U-W\|.
\]
\end{lemma}

\begin{proof}
Let $V(F)=\{1,\dots,n\}$ and $E(F)=\{e_1,\dots,e_m\}$, where
$e_t=i_tj_t$, $i_t<j_t$. Then
\begin{align*}
t_{x_1x_2}(F,U)&-t_{x_1x_2}(F,W)\\
&=\int\limits_{[0,1]^{n-2}} \Bigl(\prod_{ij\in E(F)}
W(x_i,x_j)-\prod_{ij\in E(F)} U(x_i,x_j)\Bigr)\,dx_3\dots\,dx_n.
\end{align*}
We can write
\[
\prod_{ij\in E(F)} W(x_i,x_j)-\prod_{ij\in E(F)}
U(x_i,x_j)=\sum_{t=1}^{m} X_t(x),
\]
where
\[
X_t(x) = \Bigl(\prod_{s=1}^{t-1} W(x_{i_s},x_{j_s})\Bigr)
\Bigl(\prod_{s=t+1}^m U(x_{i_s},x_{j_s})\Bigr)(W(x_{i_t},x_{j_t})-
U(x_{i_t},x_{j_t})).
\]
Consider the integral of a given term:
\begin{align*}
\Bigl|\int\limits_{[0,1]^{n-2}} X_t(x)\,\,dx_3\dots\,dx_n\Bigr|&\le
\int\limits_{[0,1]^n} |X_t(x)|\,\,dx_3\dots\,dx_n\\ &\le
\int\limits_{[0,1]^{n-2}} |W(x_{i_t},x_{j_t})-
U(x_{i_t},x_{j_t})|\,\,dx_3\dots\,dx_n.
\end{align*}
If $i_t=1$, then this integral is just $\|W-U\|_{x_1}$; if $i_t=2$,
then it is $\|W-U\|_{x_1}$; if $i_t\ge 3$, then it is $\|W-U\|$.
(Note that $i_t=1,j_t=2$ does not occur by hypothesis.) The first
possibility occurs $d_1$ times, the second $d_2$ times. This proves
the Lemma.
\end{proof}

\noindent{\bf Remark.} In \cite{LSz} a version of this lemma was
proved (for unlabeled graphs) where the $\ell_1$ norm was replaced by
the smaller ``rectangle norm''. Such a sharper version could be
proved here as well (but we don't need it).

\medskip

Applying lemma \ref{W-U} to all simple $2$-labeled graphs occurring
in a quantum graph, we get

\begin{corollary}\label{W-U-G}
Let $U,W:~[0,1]^2\to[0,1]$ be two symmetric functions and let $g$ be
any simple $2$-labeled quantum graph. Then there exists a constant
$c=c(g)$ depending only on $g$ such that for every $x,y\in[0,1]$,
\[
|t_{xy}(g,U)-t_{xy}(g,W)|\le c_g \max(\|U-W\|_x, \|U-W\|_y,\|U-W\|).
\]
\end{corollary}

Now we return to the proof of of the Theorem. Let $g$ be a simple
2-labeled quantum graph and assume that $g\equiv 0\pmod t$. Then in
particular
\[
t(gg,W)=\int\limits_0^1\int\limits_0^1 t_{xy}(g,W)^2\,dx\,dy =0,
\]
and hence
\begin{equation}\label{TXY}
t_{xy}(g,W) =0
\end{equation}
for almost all $x,y\in[0,1]$. Let $C\subset[0,1]$ denote the set
where this does not hold.

Next we show that
\begin{equation}\label{TXX}
t_{xx}(g,W) =0
\end{equation}
for almost all $x\in[0,1]$.

Suppose (\ref{TXX}) is false; then there is an $\eps>0$ and a set
$A\subseteq[0,1]$ with $\lambda(A)=\eps$ such that (say)
$t_{xx}(g,W)>\eps$ for all $x\in A$. Let $U$ be a continuous function
such that
\[
\|U-W\|<\frac{\eps^2}{9c_g}.
\]
From (\ref{L1X}) it follows that the set
\[
B=\{x\in[0,1]:~\|U-W\|_x>\frac{\eps}{3c_g}\}
\]
has measure less than $\eps/3$. For $x\in[0,1]$, let
\[
C_x=\{y\in[0,1]:~(x,y)\in C\},
\]
and let $D$ be the set of points $x\in[0,1]$ for which the set $C_x$
does not have measure 0. Clearly, $D$ has measure 0. Hence
$A\setminus B\setminus D$ has positive measure.

Let $x$ be a point with density 1 of the set $A\setminus B\setminus
D$. Choose any sequence $y_n\in A\setminus B\setminus C_x$ such that
$y_n\to x$ (such a sequence exists since $C_x$ has measure 0 by the
choice of $x$). Then we have by Corollary \ref{W-U-G}
\begin{equation}\label{XYN}
|t_{xy_n}(g,U)-t_{xy_n}(g,W)|<c_g
\max(\frac{\eps}{3c_g},\frac{\eps}{3c_g},\frac{\eps^2}{9c_g})=\frac{\eps}{3},
\end{equation}
since $x,y_n\notin B$, and similarly,
\begin{equation}\label{XX}
|t_{xx}(g,U)-t_{xx}(g,W)|<\frac{\eps}{3}.
\end{equation}
Here $t_{xy_n}(g,W)=0$ since $y_n\notin C_x$, and $t_{xx}(g,W)>\eps$,
since $x\in A$. Thus (\ref{XYN}) and (\ref{XX}) imply that
\[
|t_{xy_n}(g,U)- t_{xx}(g,U)|>\eps-2\frac{\eps}{3}>\frac{\eps}{3},
\]
which is a contradiction, since $U$ is continuous, and therefore
\[
t_{xy_n}(g,U)\to t_{xx}(g,U).
\]
This contradiction proves (\ref{TXX}).

From here, the proof of Theorem \ref{W-PROP} is easy. Trivially
\[
t_x(g',W)=t_{xx}(G,W),
\]
and so (\ref{TXX}) implies that for every $1$-labeled quantum graph
$h$
\[
t(g'h,W)=\int\limits_0^1 t_x(g'h,W)\,dx = \int\limits_0^1
t_x(g',W)t_x(h,W)\,dx=0.
\]

To prove the second assertion of the theorem, it suffices to note
that if $t$ had a contractor, then it would have a representation in
the form of $t=t(.,H)$ with some finite weighted graph $H$ by Theorem
\ref{CONTRACT}. In other words, we would have a stepfunction $W'$
such that $t(F,W)=t(F,W')$ for every finite graph $F$. By the results
of \cite{BCL}, this implies that $W$ is a stepfunction (up to set of
measure $0$).


\begin{thebibliography}{99}

\bibitem{BCL}
C.~Borgs, J.~Chayes and L.~Lov\'asz: Unique limits of dense graph
sequences (manuscript)

\bibitem{FLS}
M.~Freedman, L.~Lov\'asz, A.~Schrijver: Reflection positivity, rank
connectivity, and homomorphism of graphs (MSR Tech Report \#
MSR-TR-2004-41)

\url{ftp://ftp.research.microsoft.com/pub/tr/TR-2004-41.pdf}

\bibitem{LL}
L.~Lov\'asz: The rank of connection matrices and the dimension of
graph algebras, {\it Eur.\ J.~Combinatorics} (to appear).

\bibitem{LSz}
L.~Lov\'asz and B.~Szegedy: Limits of dense graph sequences (MSR Tech
Report \# MSR-TR-2004-79)

\url{ftp://ftp.research.microsoft.com/pub/tr/TR-2004-79.pdf}

\bibitem{LTS}
L.~Lov\'asz and V.T.~S\'os: Generalized quasirandom graphs
(manuscript).

\bibitem{FLW}
M.~Freedman, L.~Lov\'asz, D.~Welsh (unpublished)
\end{thebibliography}
\end{document}